\documentclass[11pt, twoside]{article}
\usepackage{latexsym}
\usepackage{amsmath}
\usepackage{amssymb}
\usepackage[all]{xy}
\usepackage{amsfonts}
\usepackage{verbatim}
\usepackage{amsthm}
\usepackage{mathrsfs}\usepackage{array}
\usepackage{epsfig}
\usepackage{xy}
\usepackage{tensor}
\usepackage{array}
\usepackage{stmaryrd}
\usepackage{graphicx,color}
\usepackage{xcolor}
\usepackage[colorlinks=true,linkcolor=blue,citecolor=blue]{hyperref}
\usepackage{tikz}
\usetikzlibrary{arrows,calc}
\usepackage{etex}
\usepackage{mathdots}
\usepackage{float}
\usepackage{graphics}
\usepackage{pdflscape}
\usepackage{extarrows}
\usepackage{anysize,hyperref}
\input xypic
\xyoption{all}
\usepackage{bm}
\usepackage[perpage,symbol]{footmisc}
\usepackage{setspace}
\SelectTips{cm}{10}

\def\dr{\ar@{->}[r]}
\def\Ker{\mbox{\rm Ker}\,}

\begin{document}
\baselineskip=15pt
\title{\Large{\bf  Several results on exact sequences in categories of modules over trusses
\footnotetext{$^\ast$Corresponding author.  ~Yongduo Wang is supported by the National Natural Science Foundation of China (Grant No. 10861143). ~Jian He is supported by Youth Science and Technology Foundation of Gansu Provincial (Grant No. 23JRRA825). ~Dejun Wu is supported by the National Natural Science Foundation of China (Grant No. 12261056).} }}
\medskip
\author{Yongduo Wang$^\ast$,~Dengke Jia,~Jian He and Dejun Wu}

\date{}

\maketitle
\def\blue{\color{blue}}
\def\red{\color{red}}

\newtheorem{theorem}{Theorem}[section]
\newtheorem{lemma}[theorem]{Lemma}
\newtheorem{corollary}[theorem]{Corollary}
\newtheorem{proposition}[theorem]{Proposition}
\newtheorem{conjecture}{Conjecture}
\theoremstyle{definition}
\newtheorem{definition}[theorem]{Definition}
\newtheorem{question}[theorem]{Question}
\newtheorem{notation}[theorem]{Notation}
\newtheorem{remark}[theorem]{Remark}
\newtheorem{remark*}[]{Remark}
\newtheorem{example}[theorem]{Example}
\newtheorem{example*}[]{Example}

\newtheorem{construction}[theorem]{Construction}
\newtheorem{construction*}[]{Construction}

\newtheorem{assumption}[theorem]{Assumption}
\newtheorem{assumption*}[]{Assumption}

\baselineskip=17pt
\parindent=0.5cm

\begin{abstract}
\noindent Categorical aspects of the theory of modules over trusses were studied in recent years. The snake lemma and the nine lemma in  categories of modules over trusses are formulated in this paper.\\[0.3cm]
\textbf{Keywords:} truss; snake lemma;  nine lemma\\[0.1cm]
\textbf{2020 Mathematics Subject Classification:} 18G80; 18E10
\medskip
\end{abstract}
\medskip

\pagestyle{myheadings}
\markboth{\rightline {\scriptsize Yongduo Wang$^\ast$\hspace{2mm},~Dengke Jia,~Jian He and Dejun Wu}}
         {\leftline{\scriptsize  Several results on exact sequences in categories of modules over trusses}}

\section{Introduction}
 In the 1920s, H. Pr$\ddot{u}$fer \cite{hp}, R. Baer \cite{rb} and A. K. Su$\check{s}$kevi$\check{c}$ \cite{aks} introduced a novel
algebraic structure called heap. A heap is a set with a ternary operation limited by
certain conditions, called associativity and Ma$l^{^{,}}$ce$v^{^{,}}$s identities. It
turns out that there is a deep connection between groups and heaps: heaps correspond
to free transitive actions of groups on sets, which allows us to depart from the choice
of a neutral element and focus on the set as an affine version of the group itself.

   Nearly 100 years later, trusses were introduced in \cite{tb} as structures describing two different distributive laws: the well-known ring distributivity and the one coming from the recently introduced braces, which are gaining popularity due to their role in the study of the set-theoretic solutions of the Yang-Baxter equation. The brace distributive law appeared earlier in the context of quasi-rings of radical rings; see \cite{ch}. It turns out that both these structures, rings and braces, can be described elegantly by switching the group structure to a heap structure. This leads to the definition of a truss, which is a set $T$ with a ternary operation $[-, -, -]$ and a binary multiplication $\cdot$ satisfying some conditions, the crucial one being the generalisation of ring and brace distributivity: $a\cdot [b, c, d] = [a\cdot b, a\cdot c, a\cdot d]$ and $[b, c, d]\cdot a = [b\cdot a, c\cdot a, d\cdot a]$, for all $a, b, c, d\in T$. Thanks to this, we can jointly approach brace and ring theory.

   Every truss $T$ is a congruence class of a ring $R(T)$, the universal extension of $T$ into a ring (see \cite{rtb}). Trusses, even though similar to rings, differ significantly as the category of trusses has no zero object. Having the structure of a truss, it is natural to ask: what is the theory of modules over trusses? The notion of modules over trusses was posed and basic properties of it were given by T. Brzezi$\acute{n}$ski (see \cite{tb}). In recent years, modules over trusses were studied by S. Breaz, T. Brzezi$\acute{n}$ski, B. Rybolowicz and P. Saracco from different aspects (see \cite{stbp,stbp1,tb1,tb2,tbp}). In \cite{stbp1}, they gave the concept of exact sequences in categories of modules over trusses. Here we continue to study exact sequences in categories of modules over trusses. The versions of the snake lemma and the nine lemma in  categories of modules over trusses are formulated in this paper.

\section{Preliminaries}
A heap is a set $H$ together with a ternary operation
$[---]: H\times H\times H\rightarrow H$ which is associative and satisfies the Mal'cev identities, that is,
\begin{equation*}
[[a,b,c],d,e] = [a,b,[c,d,e]] \qquad \text{and} \qquad [a,b,b] = a = [b,b,a]
\end{equation*}
for all $a,b,c,d,e\in H$.
A heap $H$ is said to be abelian if for all $a,b,c\in H$, $[a,b,c] = [c,b,a]$.

A heap morphism from $(H,[---])$ to $(H^\prime,[---])$ is a function $f:H\rightarrow H^\prime$ respecting the ternary operations, i.e., such that for all $x$, $y$, $z\in H$, $f([x,y,z])=[f(x),f(y),f(z)]$.
The category of heaps is denoted by Heap and the category of abelian heaps is denoted by Ah. A singleton set $\{\ast\}$  with the unique heap operation $[\ast,\ast,\ast]=\ast$, it is the terminal objective in the category of heaps, we denote it by $\star$. The empty set is the initial object. There is no zero objective in the category of heaps.

With every group $G$ we can associate a heap $H(G)=(G,[-,-,-])$ where $[a,b,c]=ab^{-1}c$ and every morphism of group is automatically a morphism of heaps. With every non-empty heap $H$ and for a fixed $e\in H$, we can associate a group $G(H;e)$ and the binary operation is $a\cdot b=[a,e,b]$. The inverse of $a\in G(H;e)$ is $a^{-1}=[e,a,e]$.

\begin{lemma}\label{le.se}{\rm \cite[Lemma 2.3]{tb1}}
Let $(H,[-,-,-])$ be a heap
\begin{enumerate}
\item[{\rm (1)}] If $e$, $x$, $y\in H$ are such that $[x,y,e]=e$ or $[e,x,y]=e$, then $x=y$.
\item[{\rm (2)}] For all $v$, $w$, $x$, $y$, $z\in H$
$$[v,w,[x,y,z]]=[v,[y,x,w],z].$$
\item[{\rm (3)}] For all $x$, $y$, $z\in H$
$$[x,y,[y,x,z]]=[[z,x,y],y,x]=[x,[y,z,x],y]=z.$$
In particular, in the expression $[x,y,z]=w$ any three elements determine the fourth one.
\item[{\rm (4)}] If $H$ is abelian, then, for all $x_i$, $y_i$, $z_i\in H$, $i=1,2,3$,
$$[[x_1,x_2,x_3],[y_1,y_2,y_3],[z_1,z_2,z_3]]=[[x_1,y_1,z_1],[x_2,y_2,z_2],[x_3,y_3,z_3]].$$
\end{enumerate}
\end{lemma}

A subset $S$ of a heap $H$ that is closed under the heap operation is called a sub-heap of $H$. Every non-empty sub-heap $S$ of an abelian heap $H$ defines a congruence relation $\sim_S$ on $H$:
\begin{equation*}
a\sim_S b \quad \iff \quad \exists s\in S,\ [a,b,s]\in S \quad \iff \quad \forall s\in S,\ [a,b,s]\in S.
\end{equation*}
The equivalence classes of $\sim_S$ form an abelian heap with operation induced from that in $H$. Namely,
$
[\bar a, \bar b, \bar c] = \overline{[a,b,c]}
$, where $\bar x$ denotes the class of $x$ in $H/\sim_S$ for all $x\in H$. This is known as the {\em quotient heap} and it is denoted by $H/S$.
For any $s\in S$ the class of $s$ is equal to $S$.

If $\varphi: H\rightarrow K$ is a morphism of abelian heaps, then for all $e\in {\rm Im}\varphi$ the set
\begin{equation*}
\ker_e \varphi:=  \{ a\in H\; |\; \varphi(a)=e\}
\end{equation*}
is a sub-heap of $K$. Different choices of $e$ yield isomorphic sub-heaps and the quotient heap $H/\ker_e \varphi$ does not depend on the choice of $e$. Moreover, the sub-heap relation $\sim_{\ker_e \varphi}$ is the same as the kernel relation defined by: $a\,\Ker\varphi\,b$ if and only if $\varphi(a) = \varphi(b)$. Thus we write $\ker\varphi$ for $\ker_e \varphi$ and we refer to it as the kernel of $\varphi$.

\begin{definition}{\rm \cite[Definition 3.1]{tb1}}
A truss is an algebraic system consisting of a set $T$, a ternary operation $[-,-,-]$ making $T$ into an Abelian heap, and an associative binary operation $\cdot$ which distributes over $[-,-,-]$, that is, for all $w$, $x$, $y$, $z\in T$,
\begin{equation*}
w[x,y,z]=[wx,wy,wz], \quad  [x,y,z]w=[xw,yw,zw].
\end{equation*}
A truss is said to be commutative if the binary operation $\cdot$ is commutative.
\end{definition}
A heap homomorphism between two trusses is a truss homomorphism if it respects multiplications. The category of trusses and their morphism is denoted by ${\rm Trs}$.

Let $T$ be a truss. A left $T$-module is an abelian heap $M$ together with an associative left action $\lambda_M: T\times M \rightarrow M$ of $T$ on $M$ that distributes over the heap operation. The action is denoted on elements by $t\cdot m = \lambda_M(t,m)$, with $t\in T$ and $m\in M$. Explicitly, the axioms of an action state that, for all $t,t',t''\in T$ and $m,m',m''\in M$,
\begin{subequations}
\begin{equation*}
t\cdot(t'\cdot m) = (tt')\cdot m,
\end{equation*}
\begin{equation*}
 [t,t',t'']\cdot m = [t\cdot m,t'\cdot m,t''\cdot m] ,
\end{equation*}
\begin{equation*}
 t\cdot [m,m',m''] = [t\cdot m,t\cdot m',t\cdot m''].
\end{equation*}
\end{subequations}
If $T$ is a unital  truss and the action satisfies $1\cdot m = m$, then we say that $M$ is a unital or normalised module. A sub-heap $N$ of a left $T$-module $M$ is called a submodule if it is closed under the $T$-action.

A module homomorphism is a homomorphism of heaps between two modules that also respects the actions. As it is customary in ring theory we often refer to homomorphisms of $T$-modules as to $T$-linear maps or morphisms.
The category of left $T$-modules is denoted by $T$-${\rm mod}$, that of left unital $T$-modules  by $ T_{1}$-${\rm mod}$ and the heaps of homomorphisms between modules $M$ and $N$ are denoted by ${\rm Hom}_{T}(M,N)$. The terminal heap $\star$ and initial heap $\varnothing$, with the unique possible actions, are the terminal and the initial object in $T$-${\rm mod}$. It should be noted that, since $\star\neq\varnothing$, $T$-${\rm mod}$ do not have zero object.

An element $e$ of a left $T$-module $M$ is called an absorber provided that
\begin{equation*}
t\cdot e = e, \qquad \mbox{for all $t\in T$}.
\end{equation*}
The set of all absorbers in $M$  is denoted by ${\rm Abs}(M) = \{ m\in M\mid t\cdot m = m, \forall\,t\in T\}$.

\begin{proposition}{\rm \cite[Proposition 2.6]{tbp}}
Every epimorphism of $T$-module is surjective.
\end{proposition}

\begin{proposition}{\rm \cite[Proposition 2.8]{tbp}}
Every monomorphism of $T$-module is injective.
\end{proposition}

\begin{definition}{\rm \cite[Definition 2.5]{stbp1}}
Let $M$ be a non-empty left $T$-module. For every $e \in M$, the action $\cdot_e\colon  T \times M \rightarrow M$, given by
\[t\cdot_e m = [t\cdot m,t\cdot e,e], \qquad \textrm{for all }m\in M, t \in T,\]
is called the \emph{$e$-induced action} or the \emph{$e$-induced module structure} on $M$ and denote it by $M^{(e)}$. We say that a subset $N\subseteq M$ is {\em an induced submodule} of $M$ if $N$ is a non-empty sub-heap of $M$ and $t\cdot_e n \in N$ for all $t\in T$ and $n,e\in N$.
\end{definition}
Different choices of $e$ yield isomorphic induced modules and an iteration of an induced action gives an induced action. For all $T$-module morphisms $\varphi:M\rightarrow N$, this yields an analogue of the fist isomorphism theorem for $T$-module: $M/\Ker\varphi\cong {\rm Im}\varphi$

If $R$ is a ring then we can consider its associated truss ${\rm T}(R)=({\rm H}(R,+),\cdot)$. Moreover, any $R$-module $M$ gives rise, in the same way obtain a ${\rm T}(R)$-module ${\rm T}(M)=({\rm H}(M,+),\cdot)$, whose underlying abelian heap structure is induced by the abelian group one. This assignment gives rise to a functor
$${\rm T} : R\mbox{-}{\rm mod}\longrightarrow {\rm T}(R)\mbox{-}{\rm mod}, \quad (M,+,\cdot) \longmapsto (H(M,+),\cdot), \quad f \longmapsto  f$$

Let $T$ be a truss (not necessarily unital) and let $\star$ denote the singleton $T$-module. We say that a sequence of non-empty $T$-modules $\xymatrix{M\ar[r]^{f} & N\ar[r]^{g} & P}$  is exact provided there exists $e\in {\rm Im}g$ such that ${\rm Im}f=\ker_{e}g$ as sets.

\begin{lemma}\label{lem.abs}{\rm \cite[Lemma 6.1]{tbp}}
Let $M,N,P$ be $T$-modules and $f:M \rightarrow N$ and  $g:N \rightarrow P$ be $T$-linear maps. There exist exact sequences
\begin{equation*}
\xymatrix{M\ar[r]^f & N \ar[r]^{g} & P}, \quad \xymatrix{\star \ar[r] & M^{(e)} \ar[r]^{f} & N^{(f(e))}} \quad \mbox{and}\quad  \xymatrix{N \ar[r]^{g} & P \ar[r] & \star}
\end{equation*}
if and only if
\begin{enumerate}
\item[{\rm (1)}]
$f$ is injective and
\item[{\rm (2)}]  $N/{\rm Im}f\cong P$ as $T$-modules,
\end{enumerate}
 where the module structure on $N/{\rm Im}f$ is the one for which the canonical projection $\pi:N\rightarrow N/{\rm Im}f$ is $T$-linear.
\end{lemma}

By abuse of terminology, we will say that
\[
\xymatrix{
\star \ar@{.>}[r] & M \ar@<+0.3ex>@{.>}[r]^{f} \ar@<-0.3ex>[r] & N \ar[r]^g & P \ar[r] & \star
}
\]
is a short exact sequence of $T$-modules to mean that there exists $e\in M$ such that all three sequences are exact.

\begin{lemma}{\rm \cite[Lemma 4.10]{tb1}}
Let $T$ be a truss. $M$, $N$ are left $T$-modules. Then $M\times N$ is a $T$-module with the product heap and module structures, i.e.\
\begin{enumerate}
\item[{\rm (1)}] with the heap operation defined by
$$
\left[(m_1,n_1), (m_2,n_2), (m_3,n_3)\right] = \left([m_1,m_2,m_3],[n_1,n_2,n_3]\right),
$$
for all $m_1,m_2,m_3\in M$, $n_1,n_2,n_3 \in N$;
\item[{\rm (2)}]for all $x\in T$, $m\in M$ and $n\in N$,
$$
x\cdot (m,n) = (x\cdot m, x\cdot n).
$$
\end{enumerate}
\end{lemma}

\begin{lemma}{\rm \cite[Lemma 4.13]{tb1}}
The set ${\rm Hom}_{T}(M,N)$ is a heap with the point wise heap operation.
\end{lemma}

\begin{lemma}\label{le.qa}{\rm \cite[Lemma 4.6]{tb2}}
 Let $T$ be a (unital) truss and $R$ a (unital) ring.
 \begin{enumerate}
\item[{\rm (1)}] Let $M$ be a (unital) $\mathrm{T}(R)$-module.  Then
$\mathrm{G}(M/{\rm Abs}(M); {\rm Abs}(M))$ is a (unital)  $R$-module. We denote this $R$-module by $M_\mathrm{Abs}$.
\item[{\rm (2)}]The assignment
\begin{equation*}
	\begin{aligned}
(-)_\mathrm{Abs} &: {\rm T}(R)\mbox{-}{\rm mod}\longrightarrow R\mbox{-}{\rm mod}\\
&M\longmapsto M_\mathrm{Abs}, \quad (\varphi: M\longrightarrow M')\longmapsto (\varphi_\mathrm{Abs}: \overline{m}\longmapsto \overline{\varphi(m)})
	\end{aligned}
\end{equation*}
is a functor such that, for all $R$-modules $N$, ${\rm T}(N)_\mathrm{Abs}\cong N$.
\item[{\rm (3)}] The functor $(-)_\mathrm{Abs}$ preserves monomorphisms.
\end{enumerate}
\end{lemma}

\section{Exact sequences of modules over a truss}

\setcounter{equation}{0}
\begin{proposition}\label{thm.1dim}
Let $T$ be a truss, $f:M\rightarrow N$, $g:N\rightarrow P$ be morphisms of $T$-modules. Then the sequence \[
\xymatrix{
\star \ar@{.>}[r] & M \ar@<+0.3ex>@{.>}[r]^{f} \ar@<-0.3ex>[r] & N \ar[r]^g & P \ar[r] & \star
}
\]
is short exact if and only if the sequence is exact in $N$ and $f$ is monic, $g$ is epic.
\begin{proof}
$``\Rightarrow"$ By the definition of short exact sequences of $T$-modules, there exists $e\in M$ such that three sequences
\begin{equation*}
\xymatrix{M\ar[r]^f & N \ar[r]^{g} & P}, \quad \xymatrix{\star \ar[r] & M^{(e)} \ar[r]^{f} & N^{(f(e))}} \quad \mbox{and}\quad  \xymatrix{N \ar[r]^{g} & P \ar[r] & \star}
\end{equation*}
are exact, then $f$ is monic and $g$ is epic by Lemma~\ref{lem.abs}.

$``\Leftarrow"$  Since the sequence is exact in $N$ and $g$ is epic, then $$P={\rm Im}g\cong N/{\rm Ker}g\cong N/{\rm Im}f$$
Thus, the sequence is short exact by Lemma~\ref{lem.abs} and the definition of short exact sequence of $T$-modules.
\end{proof}
\end{proposition}

\begin{proposition}
Let $R$ be a ring, $G_1$, $G_2$, $G_3$ are $R$-modules. If $$\xymatrix{0 \ar[r] & G_1\ar[r]^{f} & G_2\ar[r]^{g} & G_3\ar[r] & 0}$$ is a short exact sequence of $R$-modules, then $$\xymatrix{\star\ar@{.>}[r] & {\rm T}(G_1)\ar@<+0.3ex>@{.>}[r]^{f} \ar@<-0.3ex>[r] & {\rm T}(G_2)\ar[r]^g & {\rm T}(G_3)\ar[r] & \star}$$ is an exact sequence of $T(R)$-modules.
\end{proposition}

\begin{remark}
Every short exact sequence of $R$-modules can be regarded as a short exact sequence of $T(R)$-modules.
\end{remark}

\begin{proposition}
Let $$\xymatrix{
\star \ar@{.>}[r] & X \ar@<+0.3ex>@{.>}[r]^{f} \ar@<-0.3ex>[r] & Y \ar[r]^g & Z \ar[r] & \star
}$$
$$\xymatrix{
\star \ar@{.>}[r] & Z \ar@<+0.3ex>@{.>}[r]^{\alpha} \ar@<-0.3ex>[r] & A \ar[r]^\beta & B \ar[r] & \star
}$$
are short exact of $T$-modules, then the sequence
$$\xymatrix{
\star \ar@{.>}[r] & X \ar@<+0.3ex>@{.>}[r]^f \ar@<-0.3ex>[r] & Y \ar[r]^{\alpha g}  & A \ar[r]^\beta & B \ar[r] & \star
}$$
is exact.
\begin{proof}
We just need proof the sequence are exact in $Y$ and $A$. Since the two sequence are exact, then there exist $e\in Z$, $e^\prime\in B$ such that ${\rm Im}f={\rm ker}_{e}g$, ${\rm Im}\alpha={\rm ker}_{e^\prime}\beta$. Consider $${\rm ker}_{\alpha(e)}\alpha g=\{y\in Y|\alpha g(y)=\alpha(e)\}$$
since $\alpha$ is injective, then $g(y)=e$, it means that $y\in {\rm ker}_{e}g$. Thus, $${\rm ker}_{\alpha(e)}\alpha g={\rm ker}_{e}g={\rm Im}f$$
Moreover, $g$ is surjective, then ${\rm Im}\alpha g={\rm Im}\alpha={\rm ker}_{e^\prime}\beta$. Therefore, the sequence are exact in $Y$ and $A$.
\end{proof}
\end{proposition}

\begin{proposition}
Let $M$ be a left $T$-module and ${\rm Abs}(M)$ is not empty, then the sequence
$$\xymatrix{
\star \ar@{.>}[r] & {\rm Abs}(M) \ar@<+0.3ex>@{.>}[r]^-{f} \ar@<-0.3ex>[r] & M \ar[r]^-{g} & M/{\rm Abs}(M) \ar[r] & \star
}$$
is exact, where $f$ is inclusion homomorphism and $g$ is canonical projection.
\begin{proof}
Since ${\rm Im}f={\rm ker}_{{\rm Abs}(M)}g$, $f$ is monic and $g$ is epic, then the sequence is exact by Proposition~\ref{thm.1dim}
\end{proof}
\end{proposition}

\begin{lemma}\label{sup}Let $T$ be a truss, $f:M\rightarrow N$ be a morphism of $T$-modules. If $e$ is an absorber of M, then $f(e)$ is an absorber of $N$ and $\emph{ker}_{f(e)}f$ is a submodule of M.
\begin{proof}Since $e$ is an absorber of M and $f$ is a morphism of $T$-modules, then for all $t\in T$, $t\cdot f(e)=f(t\cdot e)=f(e)$. This means that $f(e)$ is an absorber of $N$. Since $f(t\cdot x)=t\cdot f(x)=t\cdot f(e)=f(e)$, for all $t\in T$, $x\in {\rm ker}_{f(e)}f$, $\emph{ker}_{f(e)}f$ is a submodule of M.
\end{proof}
\end{lemma}

\begin{corollary}
If $f:M\rightarrow N$ be a morphism of $T$-modules, then $f({\rm Abs}(M))\subseteq {\rm Abs}(N)$
\end{corollary}

\begin{theorem}\label{le.we} (The Snake Lemma)~~Let $T$ be a truss. Suppose that the following diagram of $T$-modules and homomorphisms is commutative and has exact rows
$$\xymatrix{ &M^\prime \ar[r]^{\varphi} \ar[d]^{f^\prime} & M \ar[r]^{\psi} \ar[d]^{f} & M^{\prime\prime} \ar[r]\ar[d]^{f^{\prime\prime}} & \star\\
	\star\ar@{.>}[r]& N^\prime \ar@<+0.3ex>@{.>}[r]^{\varphi_1} \ar@<-0.3ex>[r] & N \ar[r]^{\psi_1} & N^{\prime\prime}
}$$
Moreover, ${\rm Abs}(M^\prime)$ is not empty, then for all $e^\prime\in {\rm Abs}(M)$, there is an exact sequence
$$\xymatrix{{\rm ker}_{f^\prime(e^\prime)}f^\prime \ar[r] & {\rm ker}_{\varphi_1 f^\prime(e^\prime)}f \ar[r] & {\rm ker}_{\psi_1\varphi_1 f^\prime(e^\prime)}f^{\prime\prime}\ar[r] & N^\prime/{\rm Im}f^\prime\ar[r] & N/{\rm Im}f\ar[r] & N^{\prime\prime}/{\rm Im}f^{\prime\prime}}(1)$$
\begin{proof}
First of all, by Lemma~\ref{sup},  ${\rm ker}_{f^\prime(e^\prime)}f^\prime$, ${\rm ker}_{\varphi_1 f^\prime(e^\prime)}f$, ${\rm ker}_{\psi_1\varphi_1 f^\prime(e^\prime)}f^{\prime\prime}$ are submodules of $M^{\prime}$, $M$, $M^{\prime}$ respectively. Consider the following diagram:
$$\xymatrix{ &{\rm ker}_{f^\prime(e^\prime)}f^\prime \ar[r]^{\varphi_0} \ar[d]^{\tau^\prime} & {\rm ker}_{\varphi_1 f^\prime(e^\prime)}f \ar[r]^{\psi_0} \ar[d]^{\tau} & {\rm ker}_{\psi_1\varphi_1 f^\prime(e^\prime)}f^{\prime\prime} \ar[d]^{\tau^{\prime\prime}}\\
&M^\prime \ar[r]^{\varphi} \ar[d]^{f^\prime} & M \ar[r]^{\psi} \ar[d]^{f} & M^{\prime\prime} \ar[r]\ar[d]^{f^{\prime\prime}} & \star\\
	\star\ar@{.>}[r]& N^\prime \ar@<+0.3ex>@{.>}[r]^{\varphi_1} \ar@<-0.3ex>[r] \ar[d]^{\pi^\prime} & N \ar[r]^{\psi_1} \ar[d]^{\pi} & N^{\prime\prime} \ar[d]^{\pi^{\prime\prime}}\\
&N^\prime/{\rm Im}f^\prime\ar[r]^{\varphi_2} & N/{\rm Im}f\ar[r]^{\psi_2} & N^{\prime\prime}/{\rm Im}f^{\prime\prime}
}$$ where $\tau^\prime$, $\tau$, $\tau^{\prime\prime}$ are inclusion homomorphisms, $\pi^\prime$, $\pi$, $\pi^{\prime\prime}$ are canonical epimorphisms.
Since the rows are exact, there exist $m^{\prime\prime}\in M^{\prime\prime}$, $n^{\prime\prime}\in {\rm Im}{\psi_1}$ such that ${\rm Im}{\varphi}={\rm ker}_{m^{\prime\prime}}\psi$, ${\rm Im}{\varphi_1}={\rm ker}_{n^{\prime\prime}}\psi_1$. Since the diagram is commutative, then $f^{\prime\prime}\psi\varphi(m^\prime)=\psi_1\varphi_1 f^\prime(m^\prime)$, for all $m^\prime\in M^\prime$, namely $f^{\prime\prime}(m^{\prime\prime})=n^{\prime\prime}$. we define that
$\varphi_0:{\rm ker}_{f^\prime(e^\prime)}f^\prime\rightarrow{\rm ker}_{\varphi_1 f^\prime(e^\prime)}f, \enspace x\mapsto\varphi(x) $. Since $f\varphi(x)=\varphi_1f^\prime(x)=\varphi_1f^\prime(e^\prime)$ for all $x\in {\rm ker}_{f^\prime(e^\prime)}f^\prime$, then $\varphi(x)\in {\rm ker}_{\varphi_1 f^\prime(e^\prime)}f$. We obtaine that $\varphi_0$ is a $T$-module homomorphism as $\varphi$ is a $T$-module homomorphism. Similarly, $\psi_0=\psi\mid_{{\rm ker}_{\varphi_1 f^\prime(e^\prime)}f}$ is also a $T$-module homomorphism. Since $$\tau\varphi_0(m^\prime)=\tau(\varphi(m^\prime))=\varphi(m^\prime)=\varphi(\tau^\prime(m^\prime))=\varphi\tau(m^\prime)$$
for all $m^\prime\in {\rm ker}_{f^\prime(e^\prime)}f^\prime$, then $\tau\varphi_0=\varphi\tau^\prime$. We can get $\tau^{\prime\prime}\varphi_0=\psi\tau$ in the same way.

Define that $\varphi_2:N^\prime/{\rm Im}f^\prime\rightarrow N/{\rm Im}f, \enspace \overline{n^\prime}\mapsto\overline{\varphi_1(n^\prime)}$. If $\overline{n_1^\prime}=\overline{n_2^\prime}$, then there exists $x^\prime\in {\rm Im}f^\prime$ such that $[n_1^\prime,n_2^\prime,x^\prime]\in{\rm Im}f^\prime$. Since $\varphi_1f^\prime=f\varphi$,  $\varphi_1({\rm Im}f^\prime)\subseteq{\rm Im}f$, then $\varphi_1(x^\prime)\in {\rm Im}f^\prime$. So $\varphi_1([n_1^\prime,n_2^\prime,x^\prime])=[\varphi_1(n_1^\prime),\varphi_1(n_2^\prime),\varphi_1(x^\prime)]\in {\rm Im}f$, namely $\overline{\varphi_1(n_1^\prime)}=\overline{\varphi_1(n_2^\prime)}$. Thus, $\varphi_2$ is well-defined. Since
\begin{equation*}
	\begin{aligned}
\varphi_2([\overline{n_1^\prime},\overline{n_2^\prime},\overline{n_3^\prime}])&=\varphi_2(\overline{[n_1^\prime,n_2^\prime,n_3^\prime]})
=\overline{\varphi_1([n_1^\prime,n_2^\prime,n_3^\prime])}\\
&=\overline{[\varphi_1(n_1^\prime),\varphi_1(n_2^\prime),\varphi_1(n_3^\prime)]}=
[\overline{\varphi_1(n_1^\prime)},\overline{\varphi_1(n_2^\prime)},\overline{\varphi_1(n_3^\prime)}]\\
&=[\varphi_2(\overline{n_1^\prime}),\varphi_2(\overline{n_2^\prime}),\varphi_2(\overline{n_3^\prime})]	\\		
		\end{aligned}
\end{equation*}
and
$$\varphi_2(t\cdot\overline{n_1^\prime})=\varphi_2(\overline{t\cdot n_1^\prime})=\overline{\varphi_1(t\cdot n_1^\prime )}
=\overline{t\cdot \varphi_1(n_1^\prime)}=t\cdot\overline{\varphi_1(n_1^\prime)}=t\cdot\varphi_2(\overline{n_1^\prime})$$
for all $\overline{n_1^\prime}$, $\overline{n_2^\prime}$, $\overline{n_3^\prime}\in N^\prime/{\rm Im}f^\prime$, $t\in T$, $\varphi_2$ is a $T$-module homomorphism. In the same way, $\psi_2:N/{\rm Im}f\rightarrow N^{\prime\prime}/{\rm Im}f^{\prime\prime}, \enspace\overline{n}\mapsto\overline{\psi_1(n)}$ is also a $T$-module homomorphism. Since $\varphi_2\pi^\prime(n^\prime)=\varphi_2(\overline{n^\prime})=\overline{\varphi_1(n^\prime)}=\pi(\varphi_1(n^\prime))=\pi\varphi_1(n^\prime)$, for all $n^\prime\in N^\prime$, then $\varphi_2\pi^\prime=\pi\varphi_1$. Similarly, $\psi_2\pi=\pi^{\prime\prime}\psi_1$.

Secondly, it suffices to show $$\xymatrix{{\rm ker}_{f^\prime(e^\prime)}f^\prime \ar[r]^{\varphi_0} & {\rm ker}_{\varphi_1 f^\prime(e^\prime)}f \ar[r]^{\psi_0} & {\rm ker}_{\psi_1\varphi_1 f^\prime(e^\prime)}f^{\prime\prime}}$$ is exact. Since $\psi_0\varphi_0(m^\prime)=\psi\varphi(m^\prime)=m^{\prime\prime}$, then ${\rm Im}\varphi_0\subseteq{\rm ker}_{m^{\prime\prime}}\psi_0$. Let $s\in {\rm ker}_{m^{\prime\prime}}\psi_0$. Since the diagram commutes, $\psi(s)=\psi\tau(s)=\tau^{\prime\prime}\psi_0(s)=\psi_0(s)=m^{\prime\prime}$,  $s\in {\rm ker}_{m^{\prime\prime}}\psi={\rm Im}\varphi$. So there exists $s^\prime\in M^\prime$ such that $s=\varphi(s^\prime)$. But $s\in {\rm ker}_{m^{\prime\prime}}\psi_0\subseteq {\rm ker}_{\varphi_1 f^\prime(e^\prime)}f$, so $\varphi_1 f^\prime(e^\prime)=f(s)=f\varphi(s^\prime)=\varphi_1f^\prime(s^\prime)$. As $\varphi_1$ is injective, $f^\prime(s^\prime)=f^\prime(e^\prime)$, $s^\prime\in {\rm ker}_{f^\prime(e^\prime)}f^\prime$. Thus, $s=\varphi(s^\prime)=\varphi_0(s^\prime)\in{\rm Im}\varphi_0$, namely ${\rm Im}\varphi_0={\rm ker}_{m^{\prime\prime}}\psi_0$.
Next we will show that $$\xymatrix{N^\prime/{\rm Im}f^\prime\ar[r]^{\varphi_2} & N/{\rm Im}f\ar[r]^{\psi_2} & N^{\prime\prime}/{\rm Im}f^{\prime\prime}}$$
is exact. Since $\psi_2\varphi_2(\overline{n^\prime})=\psi_2(\overline{\varphi_1(n^\prime)})=\overline{\psi_1\varphi_1(n^\prime)}=\overline{n^{\prime\prime}}={\rm Im}f^{\prime\prime}$ for all $n^\prime\in N^\prime/{\rm Im}f^\prime$, then ${\rm Im}\varphi_2\subseteq{\rm ker}_{{\rm Im}f^{\prime\prime}}\psi_2$. Let $\overline{n}\in {\rm ker}_{{\rm Im}f^{\prime\prime}}\psi_2$, where $n\in N$. So $\psi_2(\overline{n})=\overline{\psi_1(n)}={\rm Im}f^{\prime\prime}$, this shows that $\psi_1(n)\in {\rm Im}f^{\prime\prime}$. Thus there exists $x^{\prime\prime}\in M^{\prime\prime}$ such that $\psi_1(n)=f^{\prime\prime}(x^{\prime\prime})$. As $\psi$ is surjective, there exists $x\in M$ such that $x^{\prime\prime}=\psi(x)$. So $\psi_1(n)=f^{\prime\prime}(x^{\prime\prime})=f^{\prime\prime}(\psi(x))=\psi_1f(x)$. Since $m^{\prime\prime}\in{\rm ker}_{n^{\prime\prime}}f^{\prime\prime}\subseteq M^{\prime\prime}$ and $\psi$ is surjective, there exists $m\in M$ such that $m^{\prime\prime}=\psi(m)$. Thus $f^{\prime\prime}(m^{\prime\prime})=f^{\prime\prime}(\psi(m))=\psi_1f(m)=n^{\prime\prime}$, this means that $f(m)\in{\rm ker}_{n^{\prime\prime}}\psi_1={\rm Im}\varphi_1$. Let $z=[f(m),f(x),n]\in N$, then
\begin{equation*}
\begin{aligned}
\mathcal[n,z,f(m)]=&[n,[f(m),f(x),n],f(m)]\\
         \overset{ Lemma~\ref{le.se}}{=}&[[n,n,f(x)],f(m),f(m)]\\
          =&[f(x),f(m),f(m)]=f(x).\\
\end{aligned}
\end{equation*}
This implies that $\overline{n}=\overline{z}$. Since $\psi_1(z)=[\psi_1(f(m)),\psi_1(f(x)),\psi_1(n)]=\psi_1(f(m))=n^{\prime\prime}$, then $z\in {\rm ker}_{n^{\prime\prime}}\psi_1={\rm Im}\varphi_1$. So there exists $z^\prime\in N^\prime$ such that $z=\varphi_1(z^\prime)$. Since $\varphi_2\pi^\prime(z^\prime)=\varphi_2(\overline{z^\prime})=\overline{\varphi_1(z^\prime)}=\overline{z}=\overline{n}\in {\rm Im}\varphi_2$, then
${\rm ker}_{{\rm Im}f^{\prime\prime}}\psi_2\subseteq{\rm Im}\varphi_2$. Therefore, the sequence is exact.

Thirdly, construct the connection homomorphism $$\delta:{\rm ker}_{n^{\prime\prime}}f^{\prime\prime}\rightarrow N^\prime/{\rm Im}f^\prime.$$
Let $a^{\prime\prime}\in {\rm ker}_{n^{\prime\prime}}f^{\prime\prime}\subseteq M^{\prime\prime}$. Since $\psi$ is surjective, there exists $a\in M$ such that $a^{\prime\prime}=\psi(a)$. As $\psi_1f(a)=f^{\prime\prime}\psi(a)=f^{\prime\prime}(a^{\prime\prime})=n^{\prime\prime}$, $f(a)\in{\rm ker}_{n^{\prime\prime}}\psi_1={\rm Im}\varphi_1$. Since $\varphi_1$ is injective, then there exists unique $a^\prime\in N^\prime$ such that $f(a)=\varphi_1(a^\prime)$. Let $\delta(a^{\prime\prime})=\pi^\prime(a^\prime)=\overline{a^\prime}$. If $a^{\prime\prime}=b^{\prime\prime}$, then there exist $b\in M$, $b^\prime\in N^\prime$ such that $\varphi_1(b^\prime)=f(b)$ and $a^{\prime\prime}=\psi(a)=\psi(b)$. For all $c\in {\rm ker}_{m^{\prime\prime}}\psi={\rm Im}{\varphi}$, there exists $c^\prime\in M^\prime$ such that $c=\varphi(c^\prime)$. Since $\psi([a,b,c])=[\psi(a),\psi(b),\psi(c)]=\psi(c)=m^{\prime\prime}$,  $[a,b,c]\in{\rm ker}_{m^{\prime\prime}}\psi={\rm Im}{\varphi}$. So there exists $d^\prime\in M^\prime$ such that $[a,b,c]=\varphi(d^\prime)$. Thus, $\varphi_1f^\prime(d^\prime)=f\varphi(d^\prime)=[f(a),f(b),f(c)]$. So
\begin{equation*}
\begin{aligned}
\varphi_1([b^\prime,a^\prime,f^\prime(d^\prime)])=&[\varphi_1(b^\prime),\varphi_1(a^\prime),\varphi_1f^\prime(d^\prime)]\\
=&[f(b),f(a),\varphi_1f^\prime(d^\prime)]\\
=&[f(b),f(a),[f(a),f(b),f(c)]]\\
\overset{ Lemma~\ref{le.se}}{=}&f(c)=f\varphi(c^\prime)=\varphi_1f^\prime(c^\prime)\\
\end{aligned}
\end{equation*}
Since $\varphi_1$ is injective,  $[b^\prime,a^\prime,f^\prime(d^\prime)]=f^\prime(c^\prime)$. This implies that $b^{\prime\prime}=a^{\prime\prime}$. Therefore, $\delta$ is well-defined. Obviously, $\delta$ is a $T$-module homomorphism.

Finally, we will show that the sequence $(1)$ is exact in ${\rm ker}_{n^{\prime\prime}}f^{\prime\prime}$. Let $t\in{\rm ker}_{\varphi_1 f^\prime(e^\prime)}f$, then $\delta\psi_0(t)=\delta(\psi(t))=\overline{t^\prime}$, where $f(t)=\varphi_1(t^\prime)=\varphi_1f^\prime(e^\prime)$. Since $\varphi_1$ is injective, then $t^\prime=f^\prime(e^\prime)$. So $\delta\psi_0(t)=\overline{t^\prime}={\rm Im}f^\prime$ implies that ${\rm Im}\psi_0\subseteq{\rm ker}_{{\rm Im}f^\prime}\delta$. Let $s^{\prime\prime}\in{\rm ker}_{{\rm Im}f^\prime}\delta$, then
$\delta(s^{\prime\prime})={\rm Im}f^\prime=\overline{s^\prime}$, where $s^{\prime\prime}=\psi(s)$, $f(s)=\varphi_1(s^\prime)$. Thus,
$s^\prime\in{\rm Im}f^\prime$. So there exists $q\in M^\prime$ such that $s^\prime=f^\prime(q)$, So $\varphi_1(s^\prime)=\varphi_1f^\prime(q)=f\varphi(q)$. Let $y=[s,\varphi(q),\varphi(e^\prime)]$, then
\begin{equation*}
\begin{aligned}
f(y)&=f([s,\varphi(q),\varphi(e^\prime)])=[f(s),f\varphi(q),f\varphi(e^\prime)]\\
&=[f(s),\varphi_1(s^\prime),\varphi_1f^\prime(e^\prime)]=\varphi_1f^\prime(e^\prime).\\
\end{aligned}
\end{equation*}
That is to say, $y\in{\rm ker}_{\varphi_1 f^\prime(e^\prime)}f$. Thus
\begin{equation*}
\begin{aligned}
\psi_0(y)&=\psi(y)=\psi([s,\varphi(q),\varphi(e^\prime)])\\
&=[\psi(s),\psi\varphi(q),\psi\varphi(e^\prime)]\\
&=[\psi(s),m^{\prime\prime},m^{\prime\prime}]=\psi(s)=s^{\prime\prime}\\
\end{aligned}
\end{equation*}
implies that ${\rm ker}_{{\rm Im}f^\prime}\delta\subseteq{\rm Im}\psi_0$. This shows that ${\rm Im}\psi_0={\rm ker}_{{\rm Im}f^\prime}\delta$.
Similarly, the sequence $$\xymatrix{{\rm ker}_{n^{\prime\prime}}f^{\prime\prime}\ar[r]^{\delta} & N^\prime/{\rm Im}f^\prime\ar[r]^{\varphi_2} & N/{\rm Im}f\ar[r]^{\psi_2} & N^{\prime\prime}/{\rm Im}f^{\prime\prime}}$$ is exact.

Thus, we get the exact sequence
$$\xymatrix{{\rm ker}_{f^\prime(e^\prime)}f^\prime \ar[r]^{\varphi_0} & {\rm ker}_{\varphi_1 f^\prime(e^\prime)}f \ar[r]^{\psi_0} & {\rm ker}_{\psi_1\varphi_1 f^\prime(e^\prime)}f^{\prime\prime}\ar[r]^{\delta} & N^\prime/{\rm Im}f^\prime\ar[r]^{\varphi_2} & N/{\rm Im}f\ar[r]^{\psi_2} & N^{\prime\prime}/{\rm Im}f^{\prime\prime}}.$$
\end{proof}
\end{theorem}
\begin{lemma} Let $T$ be a truss. Suppose that the following diagram of $T$-modules and homomorphism is commutative and has exact rows
$$\xymatrix{ &M^\prime \ar[r]^{\varphi} \ar[d]^{f^\prime} & M \ar[r]^{\psi} \ar[d]^{f} & M^{\prime\prime} \ar[r]\ar[d]^{f^{\prime\prime}} & \star\\
	\star\ar@{.>}[r]& N^\prime \ar@<+0.3ex>@{.>}[r]^{\varphi_1} \ar@<-0.3ex>[r] & N \ar[r]^{\psi_1} & N^{\prime\prime}
}$$
Moreover, ${\rm Abs}(M^\prime)$ is not empty.
\begin{enumerate}
\item[{\rm (1)}] If $f^\prime$ and $f^{\prime\prime}$ are injective, then $f$ is injective;
\item[{\rm (2)}] If $f^\prime$ and $f^{\prime\prime}$ are surjective, then $f$ is surjective;
\item[{\rm (3)}]If $f^\prime$ and $f^{\prime\prime}$ are isomorphisms, then $f$ is an isomorphism.
\end{enumerate}
\begin{proof}
(1) By Theorem~\ref{le.we}, for all $e^\prime\in M^\prime$, there exists the following exact sequence
$$\xymatrix{{\rm ker}_{f^\prime(e^\prime)}f^\prime \ar[r]^{\varphi_0} & {\rm ker}_{\varphi_1 f^\prime(e^\prime)}f \ar[r]^{\psi_0} & {\rm ker}_{\psi_1\varphi_1 f^\prime(e^\prime)}f^{\prime\prime}}.$$
Since $f^\prime$ and $f^{\prime\prime}$ are injective,  ${\rm ker}_{f^\prime(e^\prime)}f^\prime$ and ${\rm ker}_{\psi_1\varphi_1 f^\prime(e^\prime)}f^{\prime\prime}$ are single point sets. Since ${\rm Im}\varphi_0={\rm ker}_{m^{\prime\prime}}\psi_0={\rm ker}_{\varphi_1 f^\prime(e^\prime)}f$ and ${\rm Im}\varphi_0$ is a singleton set,  ${\rm ker}_{\varphi_1 f^\prime(e^\prime)}f$ is a singleton set. Thus, $f$ is injective.

(2) By Theorem~\ref{le.we}, there exists the following exact sequence
$$\xymatrix{N^\prime/{\rm Im}f^\prime\ar[r]^{\varphi_2} & N/{\rm Im}f\ar[r]^{\psi_2} & N^{\prime\prime}/{\rm Im}f^{\prime\prime}}.$$
Since $f^\prime$ and $f^{\prime\prime}$ are surjective,  $N^\prime/{\rm Im}f^\prime$ and $N^{\prime\prime}/{\rm Im}f^{\prime\prime}$ are singleton sets. As ${\rm Im}\varphi_2={\rm ker}_{{\rm Im}f^{\prime\prime}}\psi_2=N/{\rm Im}f$ and ${\rm Im}\varphi_2$ is a singleton set, then
$N/{\rm Im}f$ is a singleton set. Thus, $f$ is surjective.

(3) Following by (1) and (2)
\end{proof}
\end{lemma}

\begin{proposition} Suppose that the following diagram of $T$-modules and homomorphism is commutative and has exact rows
$$\xymatrix{ A^\prime \ar[r]^{\varphi_1} \ar[d]^{f} & A \ar[r]^{\psi_1} \ar[d]^{g} & A^{\prime\prime} \ar[r]\ar@{.>}[d]^{h} & \star\\
	B^\prime \ar[r]^{\varphi_2} & B \ar[r]^{\psi_2} & B^{\prime\prime}\ar[r] & \star~~,
}$$
then there exists a unique homomorphism $h:A^{\prime\prime}\rightarrow B^{\prime\prime}$ such that $\psi_2g=h\psi_1$. Moreover, if $f$ and $g$ are isomorphisms, then $h$ is an isomorphism.
\begin{proof}
Since rows are exact, there exist $a^{\prime\prime}\in A^{\prime\prime}$, $b^{\prime\prime}\in B^{\prime\prime}$ such that ${\rm Im}\varphi_1={\rm ker}_{a^{\prime\prime}}\psi_1$, ${\rm Im}\varphi_2={\rm ker}_{b^{\prime\prime}}\psi_2$.
Since $\psi_1$, $\psi_2$ are surjective, there exists $x\in A$ such that $x^{\prime\prime}=\psi_1(x)$ for every $x^{\prime\prime}\in A^{\prime\prime}$.

Define
\begin{equation*}
\begin{aligned}
h:A^{\prime\prime}&\longrightarrow B^{\prime\prime}\\
      x^{\prime\prime}&\longmapsto\psi_2g(x),\\
\end{aligned}
\end{equation*}
where $x^{\prime\prime}=\psi_1(x)$. If $x^{\prime\prime}=y^{\prime\prime}$, then there exists $y\in A$ such that $\psi_1(x)=\psi_1(y)$. For all $t\in {\rm ker}_{a^{\prime\prime}}\psi_1={\rm Im}\varphi_1$, there exists $t^\prime\in A^{\prime}$ such that $t=\varphi_1(t^\prime)$. Since $\psi_1([x,y,t])=[\psi_1(x),\psi_1(y),\psi_1(t)]=\psi_1(t)=a^{\prime\prime}$, $[x,y,t]\in {\rm ker}_{a^{\prime\prime}}\psi_1={\rm Im}\varphi_1$. Thus, there exists $a^\prime\in A^\prime$ such that $[x,y,t]=\varphi_1(a^\prime)$. So
\begin{equation*}
\begin{aligned}
\mathcal\psi_2g([x,y,t])&=[\psi_2g(x),\psi_2g(y),\psi_2g(t)]\\
                        &=\psi_2g\varphi_1(a^\prime)=\psi_2\varphi_2f(a^\prime)\\
                        &=b^{\prime\prime}.\\
\end{aligned}
\end{equation*}
Since $\psi_2g(t)=\psi_2g(\varphi_1(t^\prime))=\psi_2\varphi_2f(t^\prime)=b^{\prime\prime}$, $[\psi_2g(x),\psi_2g(y),b^{\prime\prime}]=b^{\prime\prime}$, that is to say, $\psi_2g(x)=\psi_2g(y)$. Therefore, $h$ is well-defined. It is easy to see that $h$ is a $T$-module homomorphism. For all $a\in A$, $h\psi_1(a)=h(\psi_1(a))=\psi_2g(a)$, then $h\psi_1=\psi_2g$. Assume that there is another homomorphism  $h^\prime:A^{\prime\prime}\rightarrow B^{\prime\prime}$ such that $h^\prime\psi_1=\psi_2g$.
Since $\psi_1$ is surjective and $h^\prime\psi_1=h\psi_1$, $h=h^\prime$.

Let $f$ and $g$ be isomorphisms. Next we will show that $h$ is an isomorphism. Since $\psi_2$ is surjective, then for all $b_{1}^{\prime\prime}\in B^{\prime\prime}$, there exists $b_1\in B$ such that $b_{1}^{\prime\prime}=\psi_2(b_1)$. As $g$ is an isomorphism, there exists a unique $a\in A$ such that $b_1=g(a)$, and then $b_{1}^{\prime\prime}=\psi_2(g(a))=h\psi_1(a)$. Therefore $h$ is surjective. Since $h(a^{\prime\prime})=\psi_2(g(a_1))$, where $a^{\prime\prime}=\psi_1(a_1)$, $a_1\in{\rm ker}_{a^{\prime\prime}}\psi_1={\rm Im}\varphi_1$, and so there exists $a_1^\prime\in A^\prime$ such that $a_1=\varphi_1(a_1^\prime)$.
Thus, $h(a^{\prime\prime})=\psi_2g(\varphi_1(a_1^\prime))=\psi_2\varphi_2f(a_1^\prime)=b^{\prime\prime}$, this means that $a^{\prime\prime}\in {\rm ker}_{b^{\prime\prime}}h$.
Let $x^{\prime\prime}\in{\rm ker}_{b{^{\prime\prime}}}h$.
Since $g(x)\in {\rm ker}_{b^{\prime\prime}}\psi_2={\rm Im}\varphi_2$, there exists $b_1^\prime\in B^\prime$ such that $g(x)=\varphi_2(b_1^\prime)$. As $f$ is an isomorphism, there exists a unique $a_2^\prime\in A^\prime$ such that $b_1^\prime=f(a_2^\prime)$.
Thus, $g(x)=\varphi_2(b_1^\prime)=\varphi_2f(a_2^\prime)=g\varphi_1(a_2^\prime)$, and so $x=\varphi_1(a_2^\prime)$. Thus, $\psi_1(x)=\psi_1(\varphi_1(a_2^\prime))=a^{\prime\prime}=x^{\prime\prime}$, and hence ${\rm ker}_b{^{\prime\prime}}h$ is a single point set. Therefore, $h$ is an isomorphism.
\end{proof}
\end{proposition}

\begin{proposition}
Suppose that the following diagram of $T$-modules and homomorphism is commutative and has exact rows
$$\xymatrix{ \star\ar@{.>}[r] &A^\prime \ar@<+0.3ex>@{.>}[r]^{\varphi_1} \ar@<-0.3ex>[r] \ar@{.>}[d]^{h} & A \ar[r]^{\psi_1} \ar[d]^{f} & A^{\prime\prime} \ar[d]^{g} \\
	\star\ar@{.>}[r] &B^\prime \ar@<+0.3ex>@{.>}[r]^{\varphi_2} \ar@<-0.3ex>[r] & B \ar[r]^{\psi_2} & B^{\prime\prime}
}$$
and there exist $a^{\prime\prime}\in A^{\prime\prime}$, $b^{\prime\prime}\in B^{\prime\prime}$ such that $g(a^{\prime\prime})=b^{\prime\prime}$,  ${\rm Im}\varphi_1={\rm ker}_{a^{\prime\prime}}\psi_1$, ${\rm Im}\varphi_2={\rm ker}_{b^{\prime\prime}}\psi_1$,
then there exists a unique homomorphism $h:A^{\prime}\rightarrow B^{\prime}$ such that $\varphi_2h=f\varphi_1$. Moreover, if $f$ and $g$ are isomorphisms, then $h$ is an isomorphism.
\begin{proof}
For all $a^\prime\in A^\prime$, $\psi_2f\varphi_{1}(a^\prime)=g\psi_1\varphi_1(a^\prime)=g(a^{\prime\prime})=b^{\prime\prime}$, then $f\varphi_{1}(a^\prime)\in {\rm ker}_{b^{\prime\prime}}\psi_2={\rm Im}\varphi_2$. Since $\varphi_2$ is injective,  there exists a unique $b^\prime$ such that $f\varphi_{1}(a^\prime)=\varphi_2(b^\prime)$. Define
\begin{equation*}
\begin{aligned}
h:A^{\prime}&\longrightarrow B^{\prime}\\
      a^{\prime}&\longmapsto b^\prime,\\
\end{aligned}
\end{equation*}
where $f\varphi(a^\prime)=\varphi_2(b^\prime)$. Obviously, $h$ is a $T$-module homomorphism and $\varphi_2h=f\varphi_1$.
If $h^\prime:A^{\prime}\rightarrow B^{\prime}$ satisfies $\varphi_2h^\prime=f\varphi_1$, then $\varphi_2h^\prime=\varphi_2h$.
Since $\varphi_2$ is injective, $h^\prime=h$.

Let $f$ and $g$ be isomorphisms. We will show that $h$ is isomorphic. Since  $\varphi_2(b^\prime)\in B$ for all $b^\prime\in B^\prime$ and $f$ is an isomorphism, there exists a unique $a\in A$ such that $\varphi_2(b^\prime)=f(a)$. Again since $g\psi_1(a)=\psi_2f(a)=\psi_2\varphi_2(b^\prime)=b^{\prime\prime}=g(a^{\prime\prime})$ and $g$ is an isomorphism,  $\psi_1(a)=a^{\prime\prime}$, namely, $a\in {\rm ker}_{a^{\prime\prime}}\psi_1={\rm Im}\varphi_1$. Thus, there exists a unique $a^\prime\in A^\prime$ such that $a=\varphi_1(a^\prime)$. So $h(a^\prime)=b^\prime$, and hence $h$ is surjective. If $h(a_1^\prime)=h(a_2^\prime)$, then $f\varphi_1(a_1^\prime)=f\varphi_1(a_2^\prime)$. Since $f$ is an isomorphism and $\varphi_1$ is injective,  $a_1^\prime=a_2^\prime$. Therefore, $h$ is an isomorphism.
\end{proof}
\end{proposition}

\begin{lemma}
(The Nine Lemma) Suppose that the following diagram of $T$-modules and homomorphisms is commutative and has exact rows:
$$\xymatrix{& \star \ar@{.>}[d] & \star \ar@{.>}[d] & \star \ar@{.>}[d] \\
 \star \ar@{.>}[r] & A^\prime \ar@<+0.3ex>@{.>}[r]^{f^\prime} \ar@<-0.3ex>[r] \ar@<+0.3ex>@{.>}[d]^{\alpha^\prime} \ar@<-0.3ex>[d] & B^\prime \ar[r]^{g^\prime} \ar@<+0.3ex>@{.>}[d]^{\beta^\prime} \ar@<-0.3ex>[d] & C^\prime \ar[r] \ar@<+0.3ex>@{.>}[d]^{\gamma^\prime} \ar@<-0.3ex>[d] & \star\\
\star \ar@{.>}[r] & A \ar@<+0.3ex>@{.>}[r]^{f} \ar@<-0.3ex>[r] \ar[d]^{\alpha} & B \ar[r]^{g} \ar[d]^{\beta} & C \ar[r]  \ar[d]^{\gamma}  & \star\\
\star \ar@{.>}[r] & A^{\prime\prime} \ar@<+0.3ex>@{.>}[r]^{f^{\prime\prime}} \ar@<-0.3ex>[r] \ar[d] & B^{\prime\prime} \ar[r]^{g^{\prime\prime}} \ar[d] & C^{\prime\prime} \ar[r] \ar[d] & \star \\
& \star & \star & \star
}$$
If the middle column is exact, then the last column is exact if and only if the first column is exact.
\begin{proof}
Since all rows are exact,  there exist $e^\prime\in C^\prime$, $e\in C$, $e^{\prime\prime}\in C^{\prime\prime}$ such that ${\rm Im}f^\prime={\rm ker}_{e^\prime}g^\prime$, ${\rm Im}f={\rm ker}_{e}g$, ${\rm Im}f^{\prime\prime}={\rm ker}_{e^{\prime\prime}}g^{\prime\prime}$. Because the diagram is commutative, $\gamma^\prime(e^\prime)=e$, $\gamma(e)=e^{\prime\prime}$.
$``\Leftarrow"$
As the first column and the middle column are exact,  there exist $s^{\prime\prime}\in A^{\prime\prime}$, $t^{\prime\prime}\in B^{\prime\prime}$ such that ${\rm Im}\alpha^\prime={\rm ker}_{s^{\prime\prime}}\alpha$,
${\rm Im}\beta^\prime={\rm ker}_{t^{\prime\prime}}\beta$ and $f^{\prime\prime}(s^{\prime\prime})=t^{\prime\prime}$.
For all $x^\prime\in {\rm ker}_{e}\gamma^\prime\subseteq C^\prime$, there exists $b^\prime\in B^\prime$ such that $x^\prime=g^\prime(b^\prime)$.
Since $g\beta^\prime(b^\prime)=\gamma^\prime g^\prime(b^\prime)=\gamma^\prime(x^\prime)=e$, $\beta^\prime(b^\prime)\in  {\rm ker}_{e}g={\rm Im}f$. Thus, there exists a unique $a\in A$ such that $\beta^\prime(b^\prime)=f(a)$, and hence $$f^{\prime\prime}\alpha(a)=\beta f(a)=\beta\beta^\prime(b^\prime)=t^{\prime\prime}=f^{\prime\prime}(s^{\prime\prime}).$$
Since $f^{\prime\prime}$ is monic, $\alpha(a)=s^{\prime\prime}$. This implies that $a\in {\rm ker}_{s^{\prime\prime}}\alpha={\rm Im}\alpha^\prime$. So there exists a unique $a^\prime\in A^\prime$ such that $a=\alpha^\prime(a^\prime)$. Since $\beta^\prime f^\prime(a^\prime)=f\alpha^\prime(a^\prime)=f(a)=\beta^\prime(b^\prime)$ and  $\beta^\prime$ is monic, $f^\prime(a^\prime)=b^\prime$. Thus, $x^\prime=g^\prime(b^\prime)=g^\prime f^\prime(a^\prime)=e^\prime$ implies that ${\rm ker}_{e}\gamma^\prime$ is a singleton set, and hence $\gamma^\prime$ is monic.

Since $g^{\prime\prime}\beta=\gamma g$ and $g^{\prime\prime}$, $\beta$, $g$ are epic,  $\gamma$ is epic.

Next we will show that the last column is exact in $C$. For all $c^\prime\in C^\prime$, there exists $b^\prime\in B^\prime$ such that $c^\prime=g^\prime(b^\prime)$, then $$\gamma\gamma^\prime(c^\prime)=\gamma\gamma^\prime g^\prime(b^\prime)=g^{\prime\prime}\beta\beta^\prime(b^\prime)=g^{\prime\prime}(t^{\prime\prime})=g^{\prime\prime}f^{\prime\prime}(s^{\prime\prime})
=e^{\prime\prime}$$ implies that ${\rm Im}\gamma^\prime\subseteq{\rm ker}_{e^{\prime\prime}}\gamma^\prime$.
Let $c\in {\rm ker}_{e^{\prime\prime}}\gamma^\prime$, then there exists $b\in B$ such that $c=g(b)$. Since $\gamma(c)=\gamma g(b)=g^{\prime\prime}\beta(b)=e^{\prime\prime}$, $\beta(b)\in {\rm ker}_{e^{\prime\prime}}g^{\prime\prime}={\rm Im}f^{\prime\prime}$. Thus, there exists a unique $a^{\prime\prime}\in A^{\prime\prime}$ such that $\beta(b)=f^{\prime\prime}(a^{\prime\prime})$. As $\alpha$ is epic, there exists $a\in A$ such that $a^{\prime\prime}=\alpha(a)$, and hence $\beta(b)=f^{\prime\prime}(a^{\prime\prime})=f^{\prime\prime}\alpha(a)=\beta f(a)$.
Since $f^{\prime\prime}(s^{\prime\prime})=t^{\prime\prime}$ and $\alpha$ is epic, there exists $s\in A$ such that $s^{\prime\prime}=\alpha(s)$ and $t^{\prime\prime}=f^{\prime\prime}(s^{\prime\prime})=f^{\prime\prime}\alpha(s)=\beta f(s)$. As $\beta([b,f(a),f(s)])=[\beta(b),\beta f(a),\beta f(s)]=\beta f(s)=t^{\prime\prime}$, $[b,f(a),f(s)]\in {\rm ker}_{t^{\prime\prime}}\beta={\rm Im}\beta^\prime$, and hence there exists a unique $b_{1}^\prime\in B^\prime$ such that $[b,f(a),f(s)]=\beta^\prime(b_{1}^\prime)$. Since
\begin{equation*}
\begin{aligned}
\gamma^\prime g^\prime(b_{1}^\prime)&=g\beta^\prime(b_{1}^\prime)=g([b,f(a),f(s)])\\
&=[g(b),gf(a),gf(s)]=[g(b),e,e]\\
&=g(b)=c,\\
\end{aligned}
\end{equation*}
${\rm ker}_{e^{\prime\prime}}\gamma^\prime\subseteq{\rm Im}\gamma^\prime$, and hence ${\rm Im}\gamma^\prime={\rm ker}_{e^{\prime\prime}}\gamma^\prime$.

$``\Rightarrow"$ Since the middle column and the last column are exact, there are $t^{\prime\prime}\in B^{\prime\prime}$, $s^{\prime\prime}\in C^{\prime\prime}$ such that ${\rm Im}\beta^\prime={\rm ker}_{t^{\prime\prime}}\beta$, ${\rm Im}\gamma^\prime={\rm ker}_{s^{\prime\prime}}\gamma$. Since $\gamma\gamma^\prime(e^\prime)=e^{\prime\prime}$ and $\gamma\gamma^\prime(c^\prime)=s^{\prime\prime}$ for every $c^\prime\in C^\prime$, $s^{\prime\prime}=e^{\prime\prime}$. As the diagram is commutative, $g^{\prime\prime}(t^{\prime\prime})=e^{\prime\prime}$, this means that $t^{\prime\prime}\in {\rm ker}_{e^{\prime\prime}}g^{\prime\prime}$. Thus, there is a unique $x^{\prime\prime}\in A^{\prime\prime}$ such that $t^{\prime\prime}=f^{\prime\prime}(x^{\prime\prime})$. Since $f^{\prime\prime}$ is monic and $f^{\prime\prime}\alpha\alpha^\prime(a^\prime)=\beta\beta^\prime f^\prime(a^\prime)=t^{\prime\prime}$ for all $a^\prime\in A^\prime$, $x^{\prime\prime}=\alpha\alpha^\prime(a^\prime)$.

Let $a^{\prime\prime}\in A^{\prime\prime}$. Since $\beta$ is epic, there exists $b\in B$ such that $f^{\prime\prime}(a^{\prime\prime})=\beta(b)$. Since $$\gamma g(b)=g^{\prime\prime}\beta(b)=g^{\prime\prime}f^{\prime\prime}(a^{\prime\prime})=e^{\prime\prime},$$ $g(b)\in {\rm ker}_{e^{\prime\prime}}\gamma$, and hence there is a unique $c^\prime\in C^\prime$ such that $g(b)=\gamma^\prime(c^\prime)$. Since $g^{\prime}$ is epic, there is  $b^\prime\in B^\prime$ such that $c^\prime=g^\prime(b^\prime)$, and hence $$g\beta^\prime(b^\prime)=\gamma^\prime g^\prime(b^\prime)=\gamma^\prime(c^\prime)=g(b).$$
Since $e\in C$ and $g$ is epic, there exists $b_1\in B$ such that $e=g(b_1)$, that is to say, $b_1\in {\rm ker}_{e}g={\rm Im}f$. Thus, there is a unique $a_1\in A$ such that $b_1=f(a_1)$. Therefore $f^{\prime\prime}\alpha(a_1)=\beta f(a_1)=\beta(b_1)$. Since $g([\beta^\prime(b^\prime),b,b_1])=[g\beta^\prime(b^\prime),g(b),g(b_1)]=g(b_1)=e$,  $[\beta^\prime(b^\prime),b,b_1]\in {\rm ker}_{e}g={\rm Im}f$. Thus, there exists a unique $a\in A$ such that $[\beta^\prime(b^\prime),b,b_1]=f(a)$. Therefore
\begin{equation*}
\begin{aligned}
f^{\prime\prime}\alpha(a)=\beta f(a)&=[\beta\beta^\prime(b^\prime),\beta(b),\beta(b_1)]\\
&=[t^{\prime\prime},f^{\prime\prime}(a^{\prime\prime}),f^{\prime\prime}\alpha(a_1)]\\
&=[f^{\prime\prime}(x^{\prime\prime}),f^{\prime\prime}(a^{\prime\prime}),f^{\prime\prime}\alpha(a_1)]\\
&=f^{\prime\prime}([x^{\prime\prime},a^{\prime\prime},\alpha(a_1)]).\\
\end{aligned}
\end{equation*}
Since $f^{\prime\prime}$ is monic, $\alpha(a)=[x^{\prime\prime},a^{\prime\prime},\alpha(a_1)]$, and hence $a^{\prime\prime}=\alpha([a_1,a,\alpha^{\prime}(a^{\prime})])$, this means that $\alpha$ is epic.

Since $g^\prime$ is epic,  there exists $b_{1}^\prime\in B^\prime$ such that $e^\prime=g^\prime(b_{1}^\prime)$, and hence $b_{1}^\prime\in {\rm ker}_{e^{\prime}}g^\prime={\rm Im}f^\prime$. Thus, there exists a unique $a^\prime\in A^\prime$ such that $b_{1}^\prime=f^\prime(a^\prime)$. Let $y^{\prime}\in{\rm ker}_{\alpha^{\prime}(a^\prime)}\alpha$. Since $\beta^\prime f^\prime(y^\prime)=f\alpha^\prime(y^\prime)=f\alpha^\prime(a^\prime)$, $f^\prime(y^\prime)\in {\rm ker}_{f\alpha^{\prime}(a^\prime)}\beta^\prime$. Since ${\rm ker}_{f\alpha^{\prime}(a^\prime)}\beta^\prime$ is a singleton set and $f^\prime(a^\prime)\in {\rm ker}_{f\alpha^{\prime}(a^\prime)}\beta^\prime$, $f^\prime(y^\prime)=f^\prime(a^\prime)$. Since $f^\prime$ is monic, $y^\prime=a^\prime$, that is to say, ${\rm ker}_{\alpha^{\prime}(a^\prime)}\alpha$ is a singleton set. Thus, $\alpha^\prime$ is monic.

Since $\alpha\alpha^\prime(a^\prime)=x^{\prime\prime}$ for all $a^\prime\in A^\prime$, ${\rm Im}\alpha^\prime\subseteq {\rm ker}_{x^{\prime\prime}}\alpha$. Let $x\in {\rm ker}_{x^{\prime\prime}}\alpha$. Since $\beta f(x)=f^{\prime\prime}\alpha(x)=f^{\prime\prime}(x^{\prime\prime})=t^{\prime\prime}$,  $f(x)\in {\rm ker}_{t^{\prime\prime}}\beta={\rm Im}\beta^\prime$, and hence there exists a unique $x^\prime\in B^\prime$ such that $f(x)=\beta^\prime(x^\prime)$. Since
$$\gamma^\prime g^\prime(x^\prime)=g\beta^\prime(x^\prime)=gf(x)=e=\gamma^\prime(e^\prime)$$
and $\gamma^\prime$ is monic, $g^\prime(x^\prime)=e^\prime$. Thus, $x^\prime\in {\rm ker}_{e^\prime}g^\prime={\rm Im}f^\prime$. So there exists a unique $a^\prime\in A^\prime$ such that $x^\prime=f^\prime(a^\prime)$, and hence $f\alpha^\prime(a^\prime)=\beta^\prime f^\prime(a^\prime)=\beta^\prime(x^\prime)=f(x)$. Therefore $x=\alpha^\prime(a^\prime)$ by the injectivity of $f$. Thus, ${\rm ker}_{x^{\prime\prime}}\alpha\subseteq{\rm Im}\alpha^\prime$. Therefore, ${\rm Im}\alpha^\prime= {\rm ker}_{x^{\prime\prime}}\alpha$.

\end{proof}
\end{lemma}

\begin{lemma}Let $T$ be a truss. Suppose that the following diagram of $T$-modules and homomorphism is commutative and has exact rows
$$\xymatrix{ A \ar[r]^{\varphi_1} \ar[d]^{\alpha} & B \ar[r]^{\varphi_2} \ar[d]^{\beta} & C \ar[r]^{\varphi_3} \ar[d]^{\gamma} & D \ar[r]^{\varphi_4}\ar[d]^{\delta} & E \ar[d]^{\varepsilon}\\
	 A^\prime \ar[r]^{\psi_1} & B^\prime \ar[r]^{\psi_2} & C^\prime \ar[r]^{\psi_3} & D^\prime \ar[r]^{\psi_4} & E^\prime~~.
}$$
\begin{enumerate}
\item[{\rm (1)}] If $\alpha$ is surjective, $\beta$ and $\delta$ are injective, then $\gamma$ is injective.
\item[{\rm (2)}] If $\varepsilon$ is injective, $\beta$ and $\delta$ are surjective, then $\gamma$ is surjective.
\item[{\rm (3)}]If $\alpha$, $\beta$, $\delta$, $\varepsilon$ are isomorphisms, then $\gamma$ is an isomorphism.
\end{enumerate}
\begin{proof}
(1) Since the diagram commutes and rows are exact,  there exists $c\in {\rm Im}\varphi_2$ such that ${\rm Im}\varphi_1={\rm ker}_{c}\varphi_2$,  $c^\prime\in {\rm Im}\psi_2$ such that ${\rm Im}\psi_1={\rm ker}_{c^\prime}\psi_2$ and $\gamma(c)=c^\prime$,  $d\in {\rm Im}\varphi_3$ such that ${\rm Im}\varphi_2={\rm ker}_{d}\varphi_3$,  $d^\prime\in {\rm Im}\psi_3$ such that ${\rm Im}\psi_2={\rm ker}_{d^\prime}\psi_3$ and $\delta(d)=d^\prime$.

 Let $t\in {\rm ker}_{c^\prime}\gamma$. Since $c^\prime\in {\rm Im}\psi_2$, there exists $b^\prime\in B^\prime$ such that $c^\prime=\psi_2(b^\prime)$. Thus, $\delta\varphi_3(t)=\psi_3\gamma(t)=\psi_3(c^\prime)=\psi_3\psi_2(b^\prime)=d^\prime$, and hence $\varphi_3(t)\in {\rm ker}_{d^\prime}\delta$. Since $\delta$ is injective, $\varphi_3(t)=d$. Therefore $t\in {\rm ker}_d\varphi_3={\rm Im}\varphi_2$. Thus there exists $b\in B$ such that $t=\varphi_2(b)$. As $\psi_2\beta(b)=\gamma\varphi_2(b)=\gamma(t)=c^\prime$, $\beta(b)\in {\rm ker}_{c^\prime}\psi_2={\rm Im}\varphi_1$. So there exists $a^\prime\in A^\prime$ such that $\beta(b)=\psi_1(a^\prime)$. Because $\alpha$ is surjective,  there exists $a\in A$ such that $a^\prime=\alpha(a)$. Since $\beta(b)=\psi_1(a^\prime)=\psi_1\alpha(a)=\beta\varphi_1(a)$ and $\beta$ is injective, $b=\varphi_(a)$. Thus, $t=\varphi_2(b)=\varphi_2\varphi_1(a)=c$. So ${\rm ker}_{c^\prime}\gamma$ is a single point set, and hence $\gamma$ is injective.

(2) Similar to (1).

(3) Following by (1) and (2).
\end{proof}
\end{lemma}

\begin{theorem}
Let $M$, $M^\prime$, $N$, and $N^\prime$ be left $T$-modules and let $f:M\rightarrow N$ is a $T$-module homomorphism.

(1) If $g:M\rightarrow M^\prime$ is epimorphic and there exist $s\in {\rm Im}g$, $t\in {\rm Im}f$ such that ${\rm ker}_{s}g\subseteq{\rm ker}_{t}f$, then there exists a unique homomorphism $h:M^\prime\rightarrow N$ such that $$f=hg.$$
Moreover, ${\rm ker}_{t}h=g({\rm ker}_{t}f)$ and ${\rm Im}h={\rm Im}f$, so that $h$ is monic if and only if ${\rm ker}_{m^\prime}f={\rm ker}_{t}f$, where $h(m^\prime)=t$ and $h$ is epic if and only if is epic.

(2) If $g:N^\prime\rightarrow N$ is monic with ${\rm Im}f\subseteq {\rm Im}g$, then there exists a unique homomorphism $h:M\rightarrow N^\prime$ such that $$f=hg.$$
Moreover, ${\rm ker}_{h(m)}h=g({\rm ker}_{f(m)}f)$ for all $m\in M$ and ${\rm Im}f=g^{-1}({\rm Im}f)$, so that $h$ is epic if and only if ${\rm Im}f={\rm Im}g$ and $h$ is monic if and only if $f$ is monic.
\begin{proof}
(1) Since $g$ is surjective, then for all $m^\prime\in M^\prime$ exists $m\in M$ such that $m^\prime=g(m)$. Define
\begin{equation*}
\begin{aligned}
h:M^{\prime}&\longrightarrow N^{\prime}\\
      m^{\prime}&\longmapsto f(m),\\
\end{aligned}
\end{equation*}
where $m^\prime=g(m)$.

If $m^\prime=m_1^{\prime}$, then there exists $m_1\in M$ such that $g(m)=g(m_1)$. Since $s\in {\rm Im}g$, there exists $s_0\in M$ such that $s=g(s_0)$. Thus $s_0\in {\rm ker}_sg\subseteq{\rm ker}_tf$, and hence $f(s_0)=t$. As $g([m,m_1,s_0])=[g(m),g(m_1),g(s_0)]=g(s_0)=s$, $[m,m_1,s_0]\in {\rm ker}_sg\subseteq{\rm ker}_tf$. Thus, $$f([m,m_1,s_0])=[f(m),f(m_1),f(s_0)]=[f(m),f(m_1),t]=t.$$ By Lemma~\ref{le.se}, $f(m)=f(m_1)$, so $h$ is well-defined. It is easy to see that $h$ is a $T$-module homomorphism and $hg=f$.

For all $m^\prime\in {\rm ker}_{t}h$, there exists $m\in M$ such that $m^\prime=g(m)$. Since $f(m)=hg(m)=h(m^\prime)=t$, $m\in {\rm ker}_{t}f$, and hence ${\rm ker}_{t}h\subseteq g({\rm ker}_{t}f)$. Let $m_0\in {\rm ker}_{t}f$. Since $hg(m_0)=f(m_0)=t$, $g(m_0)\in {\rm ker}_{t}h$. Therefore $g({\rm ker}_{t})\subseteq{\rm ker}_{t}h$. Thus, ${\rm ker}_{t}h=g({\rm ker}_{t}f)$. As $g$ is epic, ${\rm Im}f={\rm Im}h$ from $f=hg$.

If $h$ is monic, then there exists a unique $m^\prime$ such that $h(m^\prime)=t$, and so ${\rm ker}_{t}h=\{m^\prime\}$. Since $g({\rm ker}_{t}f)={\rm ker}_{t}h=\{m^\prime\}$, then $g({\rm ker}_{t}f)\subseteq{\rm ker}_{m^\prime}g$. Let $m\in{\rm ker}_{m^\prime}g$, $f(m)=hg(m)=h(m^\prime)=t$, then $m\in {\rm ker}_{t}f$, ${\rm ker}_{m^\prime}g\subseteq{\rm ker}_{t}f$. Therefore, ${\rm ker}_{m^\prime}g={\rm ker}_{t}f$, where $h(m^\prime)=t$.

If ${\rm ker}_{m^\prime}g={\rm ker}_{t}f$, where $h(m^\prime)=t$. For all $x^\prime\in {\rm ker}_{t}h\subseteq M^\prime$, there exists $x\in M$ such that $x^\prime=g(x)$. Since $h(x^\prime)=hg(x)=f(x)=t$, $x\in{\rm ker}_{t}f={\rm ker}_{m^\prime}g$, then $x^\prime=g(x)=m^\prime$. So $h$ is monic.

Since $f=hg$,  $h$ is epic if and only if is epic.

(2) Similar to (1).
\end{proof}
\end{theorem}

\begin{proposition} Suppose that the sequence
$$\xymatrix{\star\ar@{.>}[r] &M_1 \ar@<+0.3ex>@{.>}[r]^{f} \ar@<-0.3ex>[r] & M \ar[r]^{g} & M_2 \ar[r] & \star}$$ is exact and
there exists $e_2\in M_2$ such that ${\rm Im}f={\rm ker}_{e_2}g$. If ${\rm Abs}(M_1)$ is not empty and $e_2\in {\rm Abs}(M_2)$, then the following statements are equivalent:
\begin{enumerate}
\item[{\rm (1)}] There exists a homomorphism $h:M_2\rightarrow M$ such that $gh=1_{M_2}$;
\item[{\rm (2)}] There exists a homomorphism $k:M\rightarrow M_1$ such that $kf=1_{M_1}$;
\item[{\rm (3)}] There exists an isomorphism $\varphi:M_1\times M_2\rightarrow M$ such that the following diagram commutes.
\end{enumerate}
$$\xymatrix{\star\ar@{.>}[r] & M_1 \ar@<+0.3ex>@{.>}[r] \ar@<-0.3ex>[r] \ar@{=}[d] & M_1\times M_2 \ar[r] \ar[d]^{\varphi} & M_2 \ar[r] \ar@{=}[d] & \star\\
	\star\ar@{.>}[r] & M_1 \ar@<+0.3ex>@{.>}[r]^{f}  \ar@<-0.3ex>[r] & M \ar[r]^{g} & M_2 \ar[r] & \star
}$$
\begin{proof}
Define
\begin{equation*}
\begin{aligned}
\tau_1:M_1&\longrightarrow M_1\times M_2\\
      m_1&\longmapsto (m_1,e_2) \\
\end{aligned}
\end{equation*}
Since
\begin{equation*}
\begin{aligned}
\tau_1([m_1,m_1^\prime,m_1^{\prime\prime}])&=([m_1,m_1^\prime,m_1^{\prime\prime}],e_2)\\
&=[(m_1,e_2),(m_1^\prime,e_2),(m_1^{\prime\prime},e_2)]\\
&=[\tau_1(m_1),\tau_1(m_1^\prime),\tau_1(m_1^{\prime\prime})]\\
\end{aligned}
\end{equation*}
and $$\tau_1(t\cdot m_1)=(t\cdot m_1,e_2)=(t\cdot m_1,t\cdot e_2)=t\cdot (m_1,e_2)=t\cdot\tau_1(m_1)$$
for all $t\in T$, $m_1$, $m_1^\prime$, $m_1^{\prime\prime}\in M_1$, $\tau_1$ is a $T$-module homomorphism.

Define
\begin{equation*}
\begin{aligned}
\pi_2:M_1\times M_2&\longrightarrow M_2\\
      (m_1,m_2) &\longmapsto m_2 \\
\end{aligned}
\end{equation*}
Obviously, $\pi_2$ is a $T$-module homomorphism. It is easy to see that  $\tau_1$ is injective and $\pi_2$ is surjective.
As
\begin{equation*}
\begin{aligned}
{\rm ker}_{e_2}\pi_2&=\{(m_1,m_2)\in M_1\times M_2| \pi_2(m_1,m_2)=e_2,  m_1\in M_1, m_2\in M_2\}\\
&=\{(m_1,m_2)\in M_1\times M_2| m_2=e_2,  m_1\in M_1, m_2\in M_2\}\\
&=\{(m_1,e_2)\in M_1\times M_2| m_1\in M_1\}\\
&={\rm Im}\tau_1,\\
\end{aligned}
\end{equation*}
 the sequence
$$\xymatrix{\star\ar@{.>}[r] &M_1 \ar@<+0.3ex>@{.>}[r]^-{\tau_1} \ar@<-0.3ex>[r] & M_1\times M_2 \ar[r]^-{\pi_2} & M_2 \ar[r] & \star}$$ is exact.

$``(1)\Rightarrow(3)"$
Define
\begin{equation*}
\begin{aligned}
\varphi:M_1\times M_2&\longrightarrow M\\
      (m_1,m_2) &\longmapsto [f(m_1),h(e_2),h(m_2)] \\
\end{aligned}
\end{equation*}
Since
\begin{equation*}
\begin{aligned}
&\varphi([(m_1,m_2),(m_1^\prime,m_2^\prime),(m_1^{\prime\prime},m_2^{\prime\prime})])=\varphi([m_1,m_1^\prime,m_1^{\prime\prime}],[m_2,m_2^\prime,m_2^{\prime\prime}])\\
&=[f([m_1,m_1^\prime,m_1^{\prime\prime}]),h(e_2),h([m_2,m_2^\prime,m_2^{\prime\prime}])]\\
&=[[f(m_1),h(e_2),h(m_1)],[f(m_1^\prime),h(e_2),h(m_2^\prime)],[f(m_1^{\prime\prime}),h(e_2),h(m_2^{\prime\prime})]]\\
&=[\varphi(m_1,m_2),\varphi(m_1^\prime,m_2^\prime),\varphi(m_1^{\prime\prime},m_2^{\prime\prime})]\\
&\varphi(t\cdot (m_1,m_2))=\varphi(t\cdot m_1,t\cdot m_2)=[f(t\cdot m_1),h(e_2),h(t\cdot m_2)]\\
&=[f(t\cdot m_1),h(t\cdot e_2),h(t\cdot m_2)]=[t\cdot f(m_1),t\cdot h(e_2),t\cdot h(m_2)]\\
&=t\cdot [f(m_1),h(e_2),h(m_2)]=t\cdot \varphi(m_1,m_2)\\
\end{aligned}
\end{equation*}
for all $t\in T$, $(m_1,m_2)$, $(m_1^\prime,m_2^\prime)$, $(m_1^{\prime\prime},m_2^{\prime\prime})\in M_1\times M_2$,  $\varphi$ is a $T$-module homomorphism.

For all $m_1\in M_1$, $\varphi\tau_1(m_1)=\varphi(m_1,e_2)=[f(m_1),h(e_2),h(e_2)]=f(m_1)$, namely, $\varphi\tau_1=f$.

For all $(m_1,m_2)\in M_1\times M_2$
\begin{equation*}
\begin{aligned}
g\varphi(m_1,m_2)&=g([f(m_1),h(e_2),h(m_2)])\\
&=[gf(m_1),gh(e_2),gh(m_2)]\\
&=[e_2,e_2,m_2]=m_2=\pi_2(m_1,m_2),\\
\end{aligned}
\end{equation*}
then $g\varphi=\pi_2$. By Lemma 3.7, $\varphi$ is an isomorphism.

$``(2)\Rightarrow(3)"$
Define
\begin{equation*}
\begin{aligned}
\psi:M &\longrightarrow M_1\times M_2\\
      m &\longmapsto (k(m),g(m)) \\
\end{aligned}
\end{equation*}
Obviously, $\psi$ is a $T$-module homomorphism, $\psi f=\tau_1$ and $\pi_2 \psi=g$. By Lemma 3.7, $\psi$ is an isomorphism, namely, $\varphi=\psi^{-1}$.

$``(3)\Rightarrow(1), (2)"$ Let $e_1\in {\rm Abs}(M_1)$. Define
\begin{equation*}
\begin{aligned}
\tau_2:M_2 &\longrightarrow M_1\times M_2\\
      m_2 &\longmapsto (e_1,m_2) \\
\end{aligned}
\end{equation*}
\begin{equation*}
\begin{aligned}
\pi_1:M_1\times M_2 &\longrightarrow M_1\\
      (m_1,m_2) &\longmapsto m_1 \\
\end{aligned}
\end{equation*}
Obviously, $\tau_2$ and $\pi_1$ are $T$-module homomorphisms. Let $h=\varphi\tau_2$, $k=\pi_1\varphi^{-1}$, then $$gh=g\varphi\tau_2=\pi_2\tau_2=1_{M_2}$$
$$kf=\pi_1\varphi^{-1}f=\pi_1\varphi^{-1}\varphi\tau_1=\pi_1\tau_1=1_{M_1}.$$
\end{proof}
\end{proposition}

\begin{proposition}
Let $T$ be a truss, $M$, $N$ are left $T$-modules. Then the heap ${\rm Hom}_{T}(M,N)$ is a left $T$-module with the action given by
$$(t\cdot f)(m)=t\cdot f(m)$$
for all $f\in {\rm Hom}_{T}(M,N)$, $m\in M$, $t\in T$.

\end{proposition}

\begin{proposition}
Suppose that the sequence $\xymatrix{\star\ar@{.>}[r] & M \ar@<+0.3ex>@{.>}[r]^f \ar@<-0.3ex>[r] & N \ar[r]^{g} & P}$ is exact and there exists $e\in {\rm Im}g$ such that ${\rm Im}f={\rm ker}_{e}g$. If $e\in {\rm Abs}(P)$,  then
$$\xymatrix{\star\ar@{.>}[r] & {\rm Hom}_{T}(Q,M) \ar@<+0.3ex>@{.>}[r]^{h} \ar@<-0.3ex>[r] & {\rm Hom}_{T}(Q,N) \ar[r]^{l} & {\rm Hom}_{T}(Q,P)}$$
is exact for every left $T$-module $Q$, where \begin{equation*}
\begin{aligned}
h:{\rm Hom}_{T}(Q,M) &\longrightarrow {\rm Hom}_{T}(Q,N)\\
      \alpha &\longmapsto f\alpha \\
\end{aligned}
\end{equation*}
\begin{equation*}
\begin{aligned}
l:{\rm Hom}_{T}(Q,N) &\longrightarrow {\rm Hom}_{T}(Q,P)\\
      \beta &\longmapsto g\beta. \\
\end{aligned}
\end{equation*}
\begin{proof}
Define $\gamma:Q\rightarrow P, ~~q\mapsto e$. Obviously, $\gamma$ is a $T$-module homomorphism. For all $\alpha\in {\rm Hom}_{T}(Q,M)$, $q\in Q$, $lh(\alpha)(q)=gf\alpha(q)=e=\gamma(q)$, then $lh(\alpha)=\gamma$, and hence ${\rm Im}h\subseteq{\rm ker}_{\gamma}l$.

Let $\beta\in {\rm ker}_{\gamma}l$. For all $q\in Q$, $g\beta(q)=l(\beta)(q)=\gamma(q)=e$, then $\beta(q)\in {\rm ker}_{e}g={\rm Im}f$, and so there exists $m\in M$ such that $\beta(q)=f(m)$. Thus, by Theorem 3.8, we have
\begin{equation*}
\begin{aligned}
\alpha^\prime:Q &\longrightarrow M\\
      q &\longmapsto m, \\
\end{aligned}
\end{equation*}
where $\beta(q)=f(m)$ such that $\beta=f\alpha^\prime$. Since $\beta(q)=f(m)=f\alpha^\prime(q)=h(\alpha^\prime)(q)$, $\beta=h(\alpha^\prime)$, ${\rm ker}_{\gamma}l\subseteq{\rm Im}h$. So ${\rm Im}h={\rm ker}_{\gamma}l$.
Let $S$ be a $T$-module. If there exist $k$, $k^\prime :S\rightarrow {\rm Hom}_{T}(Q,M)$ such that $hk=hk^\prime$, then $fk(s)=fk^\prime(s)$ for all $s\in S$. Since $f$ is monic, $k(s)=k^\prime(s)$, this implies that $k=k^\prime$. Thus, $h$ is monic.
\end{proof}
\end{proposition}

\begin{proposition}
Let $R$ be a ring, $M$, $N$, $P$ are left ${\rm T}(R)$-modules and the sequence
$$\xymatrix{\star\ar@{.>}[r] & M \ar@<+0.3ex>@{.>}[r]^f \ar@<-0.3ex>[r] & N \ar[r]^{g} & P}$$
is exact. If there exist $e\in {\rm Abs}(P)$ such that ${\rm Im}f={\rm ker}_{e}g$ and $g({\rm Abs}(N))={\rm Abs}(P)$, then the sequence
$$\xymatrix{0 \ar[r] & M_{{\rm Abs}} \ar[r]^{f_{{\rm Abs}}}  & N_{{\rm Abs}} \ar[r]^{g_{{\rm Abs}}} & P_{{\rm Abs}} \ar[r] & 0}$$
is exact.
\begin{proof}
For all $\overline{m}\in M_{{\rm Abs}}$,
$$g_{{\rm Abs}}f_{{\rm Abs}}(\overline{m})=g_{{\rm Abs}}(\overline{f(m)})=\overline{gf(m)}=\overline{e}=\overline{0},$$
then ${\rm Im}f_{{\rm Abs}}\subseteq {\rm Ker}g_{{\rm Abs}}$.

Let $\overline{n}\in {\rm Ker}g_{{\rm Abs}}$. Since $g_{{\rm Abs}}(\overline{n})=\overline{g(n)}=\overline{0}$, then $g(n)\in {\rm Abs}(P)$.
As $g({\rm Abs}(N))={\rm Abs}(P)$, there exist $n^\prime$, $n^{\prime\prime}\in {\rm Abs}(N)$ such that $g(n)=g(n^\prime)$, $e=g(n^{\prime\prime})$.
Since $$g([n^{\prime\prime},n^\prime,n])=[g(n^{\prime\prime}),g(n^\prime),g(n)]=g(n^{\prime\prime})=e,$$
then $[n^{\prime\prime},n^\prime,n]\in {\rm ker}_{e}g={\rm Im}f$, there exists $m\in M$ such that $[n^{\prime\prime},n^\prime,n]=f(m)$. Thus, $[n,f(m),n^{\prime\prime}]=n^\prime\in {\rm Abs}(N)$, namely $\overline{n}=\overline{f(m)}=f_{{\rm Abs}}(\overline{m})$. Therefore, ${\rm Im}f_{{\rm Abs}}= {\rm Ker}g_{{\rm Abs}}$.

For all $\overline{p}\in P_{{\rm Abs}}$, $p\in P$. Since $g$ is surjective, there exists $n_1\in N$ such that $p=g(n_1)$, hence $g_{{\rm Abs}}(n_1)=\overline{g(n_1)}=\overline{p}$. Thus, $g_{{\rm Abs}}$ is surjective. $f_{{\rm Abs}}$ is injective by Lemma~\ref{le.qa}. Therefore, we get the exact sequence of $R$-modules.
\end{proof}
\end{proposition}

\textbf{Yongduo Wang}\\
Department of Applied Mathematics, Lanzhou University of Technology, 730050 Lanzhou, Gansu, P. R. China\\
E-mail: \textsf{ydwang@lut.edu.cn}\\[0.3cm]
\textbf{Dengke Jia}\\
Department of Applied Mathematics, Lanzhou University of Technology, 730050 Lanzhou, Gansu, P. R. China\\
E-mail: \textsf{1719768487@QQ.com}\\[0.3cm]
\textbf{Jian He}\\
Department of Applied Mathematics, Lanzhou University of Technology, 730050 Lanzhou, Gansu, P. R. China\\
E-mail: \textsf{jianhe30@163.com}\\[0.3cm]
\textbf{Dejun Wu}\\
Department of Applied Mathematics, Lanzhou University of Technology, 730050 Lanzhou, Gansu, P. R. China\\
E-mail: \textsf{wudj@lut.edu.cn}\\[0.3cm]

\end{document}